\documentstyle[12pt,twoside]{amsart}

\renewcommand{\c}[0]{{\mathbb C}}  

\renewcommand{\o}[0]{{\cal O}} 
\newcommand{\z}[0]{{\mathbb Z}}

\renewcommand{\r}[0]{{\mathbb R}}

\newcommand{\p}[0]{{\mathbb P}}

\newcommand{\q}[0]{{\mathbb Q}}
\newcommand{\map}[0]{\dasharrow}
\newcommand{\qtq}[1]{\quad\mbox{#1}\quad}

\newcommand{\im}[0]{\operatorname{im}}

\newcommand{\inter}[0]{\operatorname{Int}}

\newsymbol\subsetneq 2328
\newsymbol\onto 1310
\def\into{\DOTSB\lhook\joinrel\rightarrow}
\newsymbol\twoheadrightarrow 1310

\newtheorem{thm}{Theorem}[section]

\newtheorem{lem}[thm]{Lemma}
\newtheorem{cor}[thm]{Corollary}

\newtheorem{prop}[thm]{Proposition}

\newtheorem{conj}[thm]{Conjecture}

\theoremstyle{definition}
\newtheorem{defn}[thm]{Definition}

\newtheorem{say}[thm]{}

\newtheorem{rem}[thm]{Remark}          

\newtheorem{defn-lem}[thm]{Definition--Lemma}

\theoremstyle{remark}

\begin{document}
\bibliographystyle{amsplain}

\title{The Nash Conjecture for  Nonprojective Threefolds}
\author{J\'anos Koll\'ar}

\maketitle
\tableofcontents

\section{Introduction}

In his fundamental paper on the
topology of real algebraic varieties Nash formulated the 
 following bold conjecture:

\begin{conj}\cite{nash}\label{nash.conj}
For every compact differentiable 
manifold $M$  of dimension $n$ there is a   real algebraic variety $X$
such that $M$  is diffeomorphic   to  $X(\r)$ and
\begin{enumerate}
\item $X$ is smooth,
\item $X$ is  projective,  and
\item $X$ is birational to $\p^n$. 
\end{enumerate}
\end{conj}

It turns out that this fails for surfaces \cite{Comessatti14}
and recently it was proved that this also fails in dimension 3:

\begin{thm}\cite{koll-nash}\label{main.nash.thm}
 Let $X $ be a smooth, real,  projective 3-fold birational to $\p^3$.
Assume that
 $X(\r)$ is orientable.
  Then every connected component of $X(\r)$ is among the following:
 \begin{enumerate}
\item Seifert fibered,

\item connected sum of   several copies of $S^3/\z_{m_i}$
(called lens spaces),

\item  torus bundle over $S^1$ or doubly covered by a torus bundle over
$S^1$,

\item  finitely many other possible exceptions, or

\item obtained from one of the above 3--manifolds by repeatedly taking
connected sum with 
$\r\p^3$  and $S^1\times S^2$.
\end{enumerate}
\end{thm}

Therefore it is of interest to consider weaker versions of the
conjecture.  \cite{nash} and \cite{tognoli}
proved that the conjecture holds if we drop
the condition (\ref{nash.conj}.3). In dimension 3, 
\cite{ben-mar} prove that the conjecture also holds if we drop
the condition (\ref{nash.conj}.1) instead.
The related  ``topological Nash conjecture" is solved in
dimension 3 by 
\cite{akb-king, ben-mar} and in general by \cite{mikh}. 

The aim of this paper is to  consider the
Nash conjecture without the projectivity assumption
(\ref{nash.conj}.2).
It is not completely clear how this could be done.
Allowing quasi projective varieties instead of projective ones does not
help at all.  If $Y$ is quasi projective and $Y(\r)$ is compact then
there is a projective variety $\bar Y\supset Y$ such that
$\bar Y(\r)=Y(\r)$.

To get the right concept, we have to look at compact complex
manifolds which can be obtained from $\p^n$ by a sequence of smooth blow
ups and blow downs.
These manifolds are  quite close to projective varieties, thus
 it was 
quite a surpize to me
that   the Nash conjecture holds for  
them:

\begin{thm}\label{main.thm} For every compact, connected, 
differentiable  3--manifold $M$    there is a   compact 
complex manifold $X$ which can be  obtained from  $\p^3$ by
 a sequence of smooth, real
blow ups and downs such that 
 $M$  is diffeomorphic   to  $X(\r)$. 
\end{thm}

To be precise,  a  real structure on a complex manifold
 is a
pair $(Y,\tau)$ where $Y$ is a  complex
manifold and $\tau:Y\to Y$ an antiholomorphic involution.  $Y(\r)$ 
denotes the fixed point set of $\tau$.
 The  main problem of the 
theory  of these objects
is that all reasonable names have already been taken. ``Real
analytic space" is used for something else (cf.\ \cite{narasimhan}) and
``real complex manifold" sounds goofy.

 Complex manifolds which are bimeromorphic to a projective one
have been investigated by several
authors:

\begin{defn} \cite{artin, moish}
A compact complex manifold $Y$ is called
a {\it Moishezon manifold} or an {\it
Artin algebraic space} if it is bimeromorphic to a projective variety.
\end{defn}

By a result of Chow and Kodaira,  every smooth Moishezon  surface is 
projective (cf.\
\cite[IV.5]{bpv}).
 The first nonprojective examples in dimension 3 were found by Hironaka
(see \cite[App.B.3]{hartsh}). 

\begin{defn} A {\it real Moishezon manifold} or a {\it real  algebraic
space} is an algebraic space $Y$  with 
 an antiholomorphic involution $\tau:Y\to Y$. 
\end{defn}

It is not hard to see that if $(Y,\tau)$ is a real  algebraic
space then there is a projective
real  algebraic variety $(Y',\tau')$ and a conjugation
invariant bimeromorphic map $\phi:Y\map Y'$.

The proof of (\ref{main.thm}) follows the pattern
established in \cite{ben-mar}. 

It is known that every 3--manifold can be obtained from
$S^3$ by repeated surgery along  knots. Thus we are done
if we can realize every surgery by a suitable birational transformation
of real algebraic spaces. \cite{ben-mar} realized that this 
is possible only for  a restricted class of surgeries.
They called such surgeries  {\it d\'echirures} (\ref{triv.mod.2.def}).
 It turns out that
it is better to restrict to an even smaller class of
surgeries. I call these {\it topological flops}.
A surgery along a knot is a topological flop iff
its surgery coefficient  (\ref{sur.3manif.defs}) is congruent to $1/2$
modulo 1 (\ref{topflop.def-lem}).

\cite{ben-mar} gives a complete topological classification of
3--manifolds up to topological flops. Thus we can approach
(\ref{main.thm}) in the following 2 steps.

 First  prove (\ref{main.thm})  up to flops. This was done in
 \cite{ben-mar} with the exception of one class. The 
remaining  example is
constructed in section 6.

Second, in (\ref{real.flop.exist.thm})  we prove that  
 every topological flop can be  realized
by an algebraic flop.  This is actually quite delicate. It is clear that
this can not be done without performing other birational transformations
first. The key step is to make sure that these extraneous transformations
do not change the set of real points. This is  connected with some
interesting open problems in the theory of Moishezon manifolds.
A related process of blowing up to obtain a flop was
studied by \cite{fujiki}.

\section{Surgery on Three Dimensional Manifolds}

\begin{say}\label{def.tor-stor}
[Torus and solid torus]{\ }

As a general reference, see \cite[Chap.2]{Rolfsen76}.

Let $D^2$ be the unit disc $(|u|\leq 1)\subset \c_u$ and $S^1$ the unit
circle
$(|v|= 1)\subset \c_v$. A three--manifold $T$ diffeomorphic to 
$D^2\times S^1$ is called a {\it solid torus}. 
The boundary of a 
solid torus
$$
\partial T\sim (|u|=1)\times (|v|=1)\sim S^1\times S^1
$$
is a {\it torus}. 
Any simple closed curve on the torus $S^1\times S^1$ is isotopic to
one of the form
$$
C_{a,b}:=\im[(|z|=1)\to S^1\times S^1]\qtq{given by} z\mapsto (z^a,z^b)
\qtq{for $(a,b)=1$.}
$$
The {\it meridian} of  a solid torus
is any 
 curve isotopic to  $C_{\pm 1,0}$. Note that these curves are generators
of $\ker[\pi_1(\partial T)\to \pi_1(T)]$, so their isotopy class is well
defined.  The {\it longitude} of  a solid torus
is any 
 curve isotopic to  $C_{0,\pm 1}$. The longitude depends
on the choice of a diffeomorphism between the solid torus  and
$D^2\times S^1$.

The correspondence
$$
\left(
\begin{array}{cc}
a& b\\
c&d
\end{array}
\right)
\mapsto \left[(u,v)\mapsto (u^av^b,u^cv^d)\right]
$$
gives an isomorphism
$$
GL(2,\z)\cong \frac{(\mbox{diffeomorphisms of the torus})}
{(\mbox{modulo isotopy})}.
$$
Up to isotopy, the diffeomorphisms of a solid torus are given by
$$
(u,z)\mapsto (uz^m,z^{\pm 1})\qtq{or}
(u,z)\mapsto (\bar uz^m,z^{\pm 1}),\qtq{where} m\in\z.
$$
They correspond to the subgroup
$$
\left(
\begin{array}{cc}
\pm 1& *\\
0&\pm 1
\end{array}
\right)\subset GL(2,\z).
$$
\end{say}

\begin{say}[Surgery on 3--manifolds]\label{sur.3manif.defs}{\ }

Let $M$ be a 3--manifold and $L\subset M$ a {\it knot}. (That is, an
embedded copy of $S^1$.) Assume that $M$ is orientable along $L$.
 $L$ has a tubular neighborhood which is a solid torus $T$
with bundary $\partial T\sim S^1\times S^1$. 
Set $N:=M\setminus \inter T$.
We can think of $M$ as
being glued together from these two pieces $M=T\cup N$. 
Let $\phi:\partial T\to \partial T$ be any orientation
preserving diffeomorphism. We can glue together $T$ and $N$ using $\phi$
to obtain another 3--manifold $M_{\phi}:=T\cup_{\phi}N$.  The operation
that creates $M_{\phi}$ from $M$ is 
 called a {\it surgery} along
$L$. (If  $\phi$ is  orientation
reversing, we can compose $\phi$ with an orientation
reversing diffeomorphism of $T$ to see that we do not get anything new.)
If we fix  a diffeomorphism $h:T\sim D^2\times S^1$ then 
a diffeomorphism $\phi:\partial T\to \partial T$
corresponds to  an element of $SL(2,\z)$. 

Let $\Psi:T\to T$ be a
diffeomorphism
 and $\psi:=\Psi|_{\partial T}$. It is clear that
$M$ is   diffeomorphic to $M_{\psi}$  and, more
generally, 
$M_{\phi}$ is   diffeomorphic to $M_{\phi\circ \psi}$.

In terms of the matrix of $\phi$ we see that
the surgeries corresponding to
$$
\left(
\begin{array}{cc}
a& b\\
c&d
\end{array}
\right)
\qtq{and}
\pm\left(
\begin{array}{cc}
a& b+ma\\
c&d+mc
\end{array}
\right)
$$
give diffeomorphic 3--manifolds. From this we conclude that
the rational number $a/c$ alone determines $M_{\phi}$.
$\phi$ maps the meridian $C_{1,0}: z\mapsto (z,1)$
to $C_{a,c}:z\mapsto (z^a,z^c)$.  Thus 
the image of the meridian determines the surgery.

It is importat to note that this depends on the choice of the 
isomorphism $h:T\sim D^2\times S^1$. A different choice of $h$
changes the matrix of $\phi$ to
$$
\left(
\begin{array}{cc}
1& -m\\
0&1
\end{array}
\right)
\left(
\begin{array}{cc}
a& b\\
c&d
\end{array}
\right)
\left(
\begin{array}{cc}
1& m\\
0&1
\end{array}
\right)
=
\left(
\begin{array}{cc}
a-mc& *\\
c&d+mc
\end{array}
\right).
$$
Thus  $c(\phi):=c$ and $r(\phi):=a/c \mod 1$ are independent of 
the choice of $h$. $r(\phi)\in \q/\z$ is called the {\it surgery
coefficient} of $\phi$.

In some cases, for instance when $L$ is homologous to zero, there is a
more canonical choice for $h$ and   one obtains a well defined 
surgery coefficient $r=a/c\in \q$. See \cite[9.F]{Rolfsen76} for details.
\end{say}

\begin{defn}\label{triv.mod.2.def}
   A surgery $M\mapsto M_{\phi}$ is called  {\it trivial
mod 2} if
$c(\phi)\equiv 0 \mod 2$. This notion is called  {\it d\'echirure} in
\cite{ben-mar} where it was first studied in detail.
\end{defn}

\begin{lem}\label{triv.mod.2.lem}
 A surgery $M\mapsto M_{\phi}$ is  trivial
mod 2 iff one can choose $\phi$ such that the induced map
$$
\phi_* :H_1(\partial T,\z_2)\to H_1(\partial T,\z_2)
\qtq{is the identity.}
$$
In this case there is a natural isomorphism
$H_i(M,\z_2)\to H_i(M_{\phi},\z_2)$ for every $i$.
\end{lem}

Proof. Choose a basis of $H_1(\partial T,\z)$ and let
$$
\left(
\begin{array}{cc}
a& b\\
c&d
\end{array}
\right)
$$
 be the matrix of $\phi$. By assumption $c$ is even.
Since the matrix is in $SL(2,\z)$, $a$ and $d$ are both odd.
Multiplying by
$$
\left(
\begin{array}{cc}
1& 1\\
0&1
\end{array}
\right)
$$
if necessary, we may assume that $b$ is also even. Thus
$$
\left(
\begin{array}{cc}
a& b\\
c&d
\end{array}
\right)
\equiv
\left(
\begin{array}{cc}
1& 0\\
0&1
\end{array}
\right)
\mod 2.
$$
In this case 
the Mayer--Vietoris sequences of $M=T\cup N$ and  $M_{\phi}=T\cup_{\phi}
N$ are isomorphic with $\z_2$ coefficients. This proves the
second claim.
\qed
\medskip

Unfortunately, the final conculsion of (\ref{triv.mod.2.lem}) does not
characterize surgeries which are trivial mod 2.  For instance, a surgery
with coefficient $r=a/c$ on the unknot in $S^3$
results in a lens space $S^3/\z_a$, cf.\ \cite[9.G.1]{Rolfsen76}.
Thus it does not change $H_1(M ,\z_2)$ iff $a$ is odd;
the parity of $c$ does not matter.

Nonetheless, from the point of view of real algebraic gemetry,
(\ref{triv.mod.2.def}) is the natural notion, cf.\ \cite[E.7]{ben-mar}.

\medskip

One may ask which  3--manifolds
can be obtained from a given one by surgeries which are trivial mod 2.
\cite[Thm.C]{ben-mar} gives a complete classification.  The answer in
the general case is somewhat complicated, but the orientable case is
very easy to state:

\begin{thm}\cite{ben-mar}\label{ben-mar.orient.thm}
 Let $M_1,M_2$ be two compact
orientable 3--manifolds. Then $M_2$ can be obtained from $M_1$ by a
sequence of surgeries which are trivial mod 2 iff
$h_1(M_1,\z_2)= h_1(M_2,\z_2)$.
\end{thm}

In the general case, we use the following consequence of the
classification:

\begin{thm}\cite[Thm.B and Ex.C.5]{ben-mar}\label{ben-mar.B.thm}
 Any compact
 3--manifold  $M$ can be obtained   by a
sequence of surgeries which are trivial mod 2 from a 3--manifold $N$
such that either
\begin{enumerate}
\item $N$ can be obtained form $S^3$ by repeatedly blowing  up
points and simple closed curves, or
\item $N$ is a torus bundle over $S^1$ with monodromy
$\left(
\begin{array}{cc}
0& 1\\
1&1
\end{array}
\right)\in GL(2,\z)$.
\end{enumerate}
\end{thm}

\section{Flops}

In this section we consider flops of real algebraic spaces.
  For the main theorem we need only the simplest type of flop   and so
we  restrict our attention to flops of $(-1,-1)$-curves.
The general theory of flops on  complex 3--folds is discussed
in detail in \cite{KoMo98}. 

\begin{say}[Contraction of ruled surfaces]{\ }\label{ruled.contr}

A {\it ruled surface} is a smooth, compact, complex surface
$F$ which can be written as a $\c\p^1$-bundle over a curve.
We write it as $p:F\to C$ and a typical fiber is denoted by $E$.

$F_m$ denotes the ruled surface obtained as the projectivization
of the vector bundle $\o_{\p^1}+\o_{\p^1}(m)$. 

Let $X'$ be a complex 3--fold and $F\subset X'$
a ruled surface. Let $N_{F|X'}$ denote
the normal bundle of $F$. Assume that $(E\cdot N_{F|X'})=-1$.
By a theorem of \cite{nakano} this implies that $F$ can be contracted
along $E$. That is, there is a complex 3--manifold $Y$ 
containing a curve $C\subset Y$
and a proper 
morphism $f:X'\to Y$ such that
\begin{enumerate}
\item $f$ restricts to an isomorphism
$X'\setminus F\cong Y\setminus C$, and
\item  $f|_F$ is the  projection $p:F\to  C$.
\end{enumerate}
If $X'$ is proper then so is $Y$.
It is important to emphasize that even if $X'$ is projective,
$Y$ is usually not projective. The projectivity of $Y$
is a quite subtle question, see \cite[Sec.4.2]{flipflop} for some
results.
\end{say}

\begin{say}[Flops of $(-1,-1)$-curves, complex case]{\
}\label{flops.over.c}

Let $X$ be a complex 3-manifold and $\c\p^1\cong C\subset X$ a
curve. Assume that the normal bundle of $C$ in $X$
is isomorphic to $\o_C(-1)+\o_C(-1)$. Let us blow up $C$ in $X$.
We obtain $p_1:X_1\to X$ with exceptional divisor $Q\cong \c\p^1\times
\c\p^1$. Moreover, the normal bundle   $N_{Q,X_1}$ has intersection
number $-1$ with both rulings of $Q$. 

Thus $Q\subset X_1$ is symmetrical with respect to the two factors of
$Q$ and   by (\ref{ruled.contr}) it can be contracted in the other
direction to obtain $p_2:X_1\to X'$. 
The bimeromorphic map $\phi:=p_2\circ p_1^{-1}:X\map X'$
is called the {\it flop of $C$}.   Frequently
$X'$ is also called a flop of $X$.
Note that $X$ and $X'$ have a completely symmetrical role
in a flop.
\end{say}

\begin{defn} There are two different interval bundles over $S^1$.
The trivial one is obtained from $[0,2\pi]\times[-1,1]$ by identifying
the points $(0,y)\leftrightarrow(2\pi,y)$. The total space is a
{\it cylinder}. The nontrivial 
one is obtained from $[0,2\pi]\times[-1,1]$ by identifying
the points $(0,y)\leftrightarrow(2\pi,-y)$. The total space is a
{\it M\"obius band}, it is not orientable. 
Its interior is called the open M\"obius band.
The image of the curve
$[0,2\pi]\times \{0\}$ is called the {\it center} of the
cylinder or M\"obius band. The cylinder and the  M\"obius band
both retract to their center. 

Let $L$ be a  line bundle   of degree $d$ on $\p^1$   defined over $\r$.
(With a slight abuse of notation I will write $L\cong \o_{\p^1}(d)$.)
Then $L(\r)$ is an $\r^1$-bundle over $\r\p^1\sim S^1$. 
It is the trivial bundle if $d$ is even and the open M\"obius band if $d$
is odd.
\end{defn}

\begin{say}[Flops of $(-1,-1)$-curves, real case]{\ }\label{flops.over.r}

Notation as in (\ref{flops.over.c}).
Assume in addition that $X$ has a real structure and $C$ is also real
and isomorphic to $\p^1$ over $\r$, thus $C(\r)\sim S^1$.
The isomorphism $N_{C|X}\cong \o_C(-1)+\o_C(-1)$
gives that
$$
N_{C(\r)|X(\r)}\cong \o_{C(\r)}(-1)+\o_{C(\r)}(-1)
\cong B_1+B_2
$$
 where each $B_i$ is an open M\"obius
band. 

We can identify $N_{C(\r)|X(\r)}$ with a tubular neighborhood
$ T\supset C(\r)$ such that     $\partial B_1$ is a simple
closed curve $J\subset \partial T$. 

$X_1(\r)$ is obtained from $X(\r)$ by blowing up $C(\r)$ which is
contained in $B_1$. Thus the birational transform $(p_1)^{-1}_*B_1$
is isomorphic to $B_1$ and it intersects $Q(\r)$ in 
a fiber of  $p_2$. Hence the birational transform of $B_1$ in $X'$
is a disc which intersects $C'(\r)$ transversally. 
Thus $X(\r)\mapsto X'(\r)$ is a surgery along $C(\r)$
where $J$ becomes the new meridian.
\end{say}

\begin{defn-lem}\label{topflop.def-lem}
 A surgery $M\mapsto M_{\phi}$ 
along a knot $L\subset T\subset M$ 
is called a {\it
topological flop} if the following equivalent conditions are satisfed:
\begin{enumerate}
\item $r(\phi)\equiv 1/2 \mod 1$.
\item There is an embedded M\"obius band $B\subset T$ whose center is
$L$  such that $\phi$ maps the meridian of $ T$
to $\partial B$.
\end{enumerate}
\end{defn-lem}

Proof. Let  $J\subset \partial T$ be
the image of the meridian under the surgery $\phi$. If $J$ is the
boundary of an embedded M\"obius band then 
$[J]=\pm 2[L]$ in $\pi_1(T)$, thus $c(\phi)=\pm 2$ and so 
$r(\phi)\equiv 1/2 \mod 1$. Conversely,
if
$J\subset \partial T$ is  a simple closed curve 
in a homotopy class $(a,\pm 2)$ as in (\ref{def.tor-stor})
then it is the boundary of the embedded M\"obius band
$$
[0,2\pi]\times[-1,1]\to D^2\times S^1
\qtq{given by} (x,y)\mapsto (ye^{iax/2},e^{ix}).\qed
$$

Since a surgery is determined by the image of the meridian
under $\phi$, we obtain the following:

\begin{prop}\label{topflop=mb.prop}
 Let  $L\subset M$ be  a simple closed curve such that $M$ is
orientable along $L$. Then
 there is a one--to--one correspondence
between
\begin{enumerate}
\item embedded M\"obius bands $B\subset M$ with center $L$
modulo isotopy, and
\item    topological flops of $L$.\qed
\end{enumerate}
\end{prop}

\begin{defn} Let $M$ be a 3--manifold and $B\subset M$ an embedded
M\"obius band with center $L$. 
Assume that $M$ is orientable along $L$. 
By (\ref{topflop=mb.prop})
$B$ defines a unique surgery  $M\mapsto M'$.
We call $M'$ the {\it topological flop} of $M$ along $B$.

The diffeomorphism class of $M'$ is  determined by the pair
$(M,B)$ up to isotopy.
\end{defn}

The relationship between algebraic and topological flops is
summarized in the next result which is a direct consequence of
our previous discussions.

\begin{prop}\label{alg.flop=top.flop}
 Let $X$ be a 3--dimensional, smooth, real, algebraic space
and
$\p^1\cong C\subset X$ a real curve. Assume that the normal bundle
of
$C$ in $X$ is isomorphic to $\o_C(-1)+\o_C(-1)$. 
Let $B\subset X(\r)$ be a M\"obius band corresponding to
one of the  $\o_C(-1)$ summands and $X\map X'$
 the flop of $C$.

Then $X'(\r)$ is diffeomorphic to the topological flop of $X(\r)$
along
$B$.\qed
\end{prop}

The following proposition shows that the study of
topological flops is equivalent to the study of
surgeries which are
trivial mod 2.

\begin{prop} A  surgery  $M\mapsto M'$   which is
trivial mod 2 can be written as a composite of
topological flops.
\end{prop}

Proof. Let $M\mapsto M'$ correspond to the matrix
$$
A=\left(
\begin{array}{cc}
a& b\\
c&d
\end{array}
\right)\qtq{where $c\equiv 0\mod 2$.}
$$
By (\ref{H1.gen.cor}), we can write $A$ as the product
of matrices of the form
$$
\left(
\begin{array}{cc}
\pm 1& *\\
0&\pm 1
\end{array}
\right)
\qtq{and}
\left(
\begin{array}{cc}
1& 0\\
2&1
\end{array}
\right).
$$
The first type of these corresponds to a change of the 
diffeomorphism $T\sim D^2\times S^1$
and the second of these corresponds to a topological flop by
(\ref{topflop.def-lem}).\qed

\section{Algebraic approximation of M\"obius bands}

\begin{thm}\label{moeb.band.exist.thm}
 Let $X$ be a compact, smooth, real, algebraic space
which is birational to $\p^3$. Let 
$ B\subset X(\r)$
be an embedded M\"obius band. Then there is 
\begin{enumerate}
\item 
 a smooth, real, rational curve 
$\p^1\cong C\subset X$ such that $C\cdot K_X\leq 0$, and 
\item a real subline bundle $\o_C(-1)\cong M\subset N_{C|X}$
\end{enumerate}
such that 
 the M\"obius band
corresponding to  $M$  is isotopic to $B$ in $X(\r)$.
\end{thm}

Proof.  Let $L\subset B$ denote the center of $B$. 
 Choose a tubular
neighborhood  $L\subset T\subset X(\r)$ such that
$B\subset T$ and $\partial B\subset \partial T$ is a
simple closed loop. Let $J\in H_1(\partial T,\z)$ denote its class.

If $L'\subset T$ is another knot which is isotopic to $L$
by an isotopy inside $T$ then $T$ is also a tubular neighborhood of
$L'$. Assume that $L'=C(\r)$ for some algebraic curve $C\cong \p^1$
and $\o_C(-1)\cong M\subset N_{C|X}$ is a 
real subline bundle. $M$ gives  a M\"obius band in $T$
and hence a 
simple closed loop  $J(M)\in H_1(\partial T,\z)$. 
Our aim is to find $C$ and $M$ such that $J(M)=J$.
This will be done in several steps.

By \cite{bo-ku} there is a morphism $f:\p^1\to X$ such that
$L$ is isotopic to  $f(\r\p^1)$ inside $T$. Moreover, we can assume that
$f^*T_X$ is ample.  Indeed, by \cite[2.5]{bo-ku} we can assume that
$f(\r\p^1)$ passes through any number of preassigned points of
$\inter T$. Thus $f^*T_X$ is ample by \cite[II.3.10.1]{koll96} for
general choice of
$f$.

Similarly, one can choose  $C\subset
X$ such that  
the normal bundle  $N_{C|X}$ is the sum of  two line bundles of large
degree.  This means that there are lots of injections 
$\o_C(-1)\into N_{C|X}$ and we can represent many
M\"obius bands. Unfortunately,
one also needs to find a way to 
understand how the direct summands of $N_{C|X}$ twist
compared to 
a fixed longitude of $L$.
 I was unable to find a
sensible way to do it, thus we proceed in a somewhat roundabout way.

\begin{say} Let $E$ be a rank 2 complex vector bundle on $\p^1$ defined
over
$\r$.  It induces a rank 2 real vector bundle   $E(\r)\to \r\p^1$.
Let $L\subset E(\r)$ be the zero section and 
 $L\subset T\subset E(\r)$ a tubular neighborhood.
We saw in (\ref{def.tor-stor}) that there is a canonical choice for
a meridian $m\in H_1(\partial T,\z)$ (up to sign).
The real algebraic  structure makes it possible
to choose a longitude as well in some cases.

Assume that $E(\r)$ is orientable.
We can write $E\cong \o_{\p^1}(a)+\o_{\p^1}(b)$ for some $a,b$
where $a\equiv b\mod 2$ by the orientability assumption.

Assume first that $a,b$ are both even and $a\leq b$. 
Then $(0, (s^2+t^2)^{b/2})$ is a nowhere zero  section of
$E(\r)$, thus it defines a
longitude $\ell \in H_1(\partial T,\z)$ (up to sign).
$\{m,\ell\}$ is a basis of $H_1(\partial T,\z)$.

The situation is a bit more complicated if $a\leq b$ are both odd.
In this case set 
$$
H_1(\partial T,\z)^{(2)}:=\{c\in H_1(\partial T,\z) : (m\cdot c)\equiv
0\mod 2\}.
$$
$H_1(\partial T,\z)^{(2)}$ is an index 2 subgroup of $H_1(\partial
T,\z)$. If $B\subset T$ is a M\"obius band with center
$L$ then $[\partial B]\in H_1(\partial T,\z)^{(2)}$.
The $\o_{\p^1}(b)$ summand corresponds to a such a M\"obius band.
This gives a well defined element
$\ell^{(2)} \in H_1(\partial T,\z)^{(2)}$ (up to sign)
and $\{m,\ell^{(2)}\}$ is a basis of $H_1(\partial T,\z)^{(2)}$.
In this case set 
$\ell:=\frac12 \ell^{(2)}\in H_1(\partial T,\q)$.
\end{say}

\begin{lem} Notation as above. A given  $d\in H_1(\partial T,\z)$
can be represented as the boundary of a M\"obius band
corresponding to an embedding
$\o_{\p^1}(-1)\into E$ if $d$ is not divisible by 2, 
$m\cdot d=\pm 2$ and $|\ell\cdot  d|\leq \min\{a,b\}$.
\end{lem}

Proof. Assume first that $a,b$ are both even. A map
$$
j:\o_{\p^1}(-1)\stackrel{(f,g)}{\longrightarrow}
\o_{\p^1}(a)+\o_{\p^1}(b)
$$
 is given by
a pair of homogeneous polynomials
 of degrees $a+1,b+1$. The image of $j$ is a subline bundle
iff $f$ and $g$ have no common zeros.
Let $J(j)\subset \partial T$ be the curve 
corresponding to $\im j$.

The intersection points of $J(j)$ and $\ell$ correspond to the zeros of
$f$. Choose $f$ and $g$ such that they both have $c\leq 1+\min\{a,b\}$
real roots and these roots alternate in $\r\p^1$. It is then easy to see that
$|\ell\cdot J(j)|=c$. Replacing $g$ by $-g$ will change the sign.

The case when $a,b$ are both odd are similar.\qed
\medskip

Let $f:\p^1\to X$ be a morphism which is an embedding in a 
complex analytic neighborhood of $\r\p^1$. Then $f(\r\p^1)\subset X(\r)$
is a knot and the notions of tubular neighborhood, meridian and longitude
make sense as before.
Instead of a normal bundle we have a normal sheaf $N_f:=f^*T_X/T_{\p^1}$
which can be written as
$$
N_f=\o_{\p^1}(a_f)+\o_{\p^1}(b_f)+(\mbox{torsion}),
$$
where the support of the torsion is disjoint from $\r\p^1$.
Let $C\subset X$ be a smooth rational curve 
such that $N_{C|X}\cong \o_{\p^1}(a)+\o_{\p^1}(b)$.
Consider the composite
$$
f_m:\p^1\stackrel{\phi_m}{\to} \p^1\cong C
$$
where $\phi_m$ is as in (\ref{unram.map.lem}).
Then 
$$
N_{f_m}=\o_{\p^1}(ma)+\o_{\p^1}(mb)+(\mbox{torsion}).
$$
Assume that  $a,b>0$ 
and let $B\subset X(\r)$ be a M\"obius band
with center $C(\r)$. 
Then for $m\gg 1$ there is a subline bundle
$j_m:\o_{\p^1}(-1)\into N_{f_m}$
such that  the  M\"obius band corresponding to $j_m$  is isotopic to $B$.

This is essentially what we want, except that $f_m$ is not an embedding.
We claim, however, that a small perturbation of $f_m$
gives a solution of (\ref{moeb.band.exist.thm}). 

By \cite[II.3.14.3]{koll96} there is a morphism
$F:\Delta\times \p^1\to X$ such that
$F|_{\{0\}\times \p^1}=f_m$ and $F|_{\{t\}\times \p^1}$
is an embedding for $t\neq 0$.

The   normal sheaf
$N_F:=F^*T_X/T_{\Delta\times \p^1}$ is flat over
$\Delta$ and so the injection
$$
j:\o_{\p^1}(-1)\into N_{f_m}\cong N_F|_{\{0\}\times \p^1}
$$
extends to   injections 
$$
j_t:\o_{\p^1}(-1)\into  N_F|_{\{t\}\times \p^1}
\cong N_{C_t|X}
$$
since $H^1(\p^1, N_{f_m}(1))=0$. \qed

\begin{lem}\label{unram.map.lem}
  $x\mapsto x+(x^2+1)^{-m}$
extends to a morphism  $\phi_m:\p^1\to \p^1$
of degree $2m+1$. $\phi_m:\r\p^1\to \r\p^1$ is a homeomorphism
and $\phi_m$ is unramified along $\r\p^1$.\qed
\end{lem}

\section{Algebraic realization of toplogical flops}

\begin{thm}\label{real.flop.exist.thm}
 Let $X$ be a compact, smooth, real, algebraic space and
$\p^1\cong C\subset X$ a smooth, real, rational curve such that
$C\cdot K_X\leq 0$. Let $\o_C(-1)\subset N_{C|X}$ be
a real subline bundle and $C(\r)\subset B\subset X(\r)$
the corresponding M\"obius band. 
Assume that $X(\r)$ is orientable along $C(\r)$.
Then there is
\begin{enumerate}
\item a  blow up $p:X_1\to X$ such that $X_1(\r)\to X(\r)$ is a homeomorphism,
 and
\item a flop $\pi:X_1\map X_2$ of $p^{-1}_*(C)$
\end{enumerate}
such that
\begin{enumerate}\setcounter{enumi}{2}
\item $X_2$ is a compact, smooth, real, algebraic space birational to
$X$, and
\item $X_2(\r)$ is diffeomorphic to the topological flop of
$X(\r)$ along $B$. 
\end{enumerate}
\end{thm}

Proof. $X(\r)$ is orientable along $C(\r)$, thus
$w_1(X(\r))\cdot C(\r)=0$ 
where $w_1$ denotes the first Stiefel--Whitney class
(cf.\ \cite{milnor-st}).
$$
c_1(X(\c))\cdot C\equiv w_1(X(\r))\cdot C(\r)\mod 2
$$
hence 
$$
\deg N_{C|X}=-K_X\cdot C-2=c_1(X(\c))\cdot C-2\equiv 0 \mod 2.
$$
Thus $L':=N_{C|X}/L$ is a line bundle of   degree $-1+2r$ for some $r\geq
0$. By (\ref{gen.pos.lem}) there is  a smooth, real, algebraic curve
$D\subset X$ such that
\begin{enumerate}
\item $D(\r)=\emptyset$, 
\item $D\cap C$ is precisely $r$ conjugate point pairs
$P_1,\bar P_1,\dots,P_r,\bar P_r$, and
\item the tanget vector of $D$ at every point $P_i,\bar P_i$
maps to a nonzero vector in $L\subset N_{C|X}$.
\end{enumerate}
Let us blow up $D$ to obtain  $p:X_1:=B_DX\to X$
and set $C_1:=p^{-1}_*C$. There is an exact sequence
$$
0\to L\to N_{C_1|X_1}\to  L'(-\sum (P_i+\bar P_i))\to 0
$$
Since $\deg L'(-\sum (P_i+\bar P_i))=-1+2r-2r=-1$, the sequence
splits and 
$$
N_{C_1|X_1}\cong \o_{C_1}(-1)+\o_{C_1}(-1).
$$
Thus $C_1\subset X_1$ can be flopped to get $X_2$.
(\ref{real.flop.exist.thm}.4) follows from (\ref{alg.flop=top.flop}).\qed

\begin{say}[General position curves in algebraic spaces]{\ }

If $X$ is a quasi projective variety, it is easy to
construct a curve $D$ as above by taking the intersection
of suitable hypersurfaces. On an algebraic space there may not exist
any base point free linear system, so intersections of hypersurfaces
are unlikely to produce  the needed curve $D$. 
In fact, there are singular algebraic spaces with a unique singular
point $P\in X$ such that every curve in $X$ passes through $P$.
No such smooth examples are known and it is conjectured that
this can not happen for smooth algebraic spaces.
(See \cite[5.2.6]{flipflop} for the case when $\dim X\leq 3$.)

The next lemma shows the existence of the required curve $D$
in some spceial cases which are sufficient for our purposes.
Note also that the nonprojective case is needed only 
for 3--manifolds as in (\ref{ben-mar.B.thm}.2). In section 6  we give an
explicit construction for such a real algebraic space.
It is not hard to use this construction to show the existence of the 
curve $D$.
\end{say}

\begin{lem}\label{gen.pos.lem}
 Let $X$ be a  smooth, proper, real   algebraic space
of dimension $n\geq 3$ 
which is birational to $\p^n$. Let
$Z\subset X$ be a subset of codimension
at least 2, $ p_i,\bar p_i\in X\setminus Z$   conjugate pairs of  points
and
$0\neq v_i\in T_{p_i}X$  tangent vectors. Then there is a smooth, proper,
real algebraic curve
$D\subset X$ such that
\begin{enumerate}
\item $D(\r)=\emptyset$,
\item $D\cap Z=\emptyset$,
\item $p_i\in D$ and $T_{p_i}D=\c\cdot v_i$ for every $i$.
\end{enumerate}
\end{lem}

Proof.  First we find  morphisms, defined over $\c$,
$g_i:\c\p^1\to X(\c)$ whose image passes through $p_i$ and
has tangent space $\c\cdot v_i$ there. Then we prove that
if we choose the $g_i$ sufficiently general then
the union of $\im g_i$ and their conjugates is  a curve $D$ as desired.
A modification of this construction would give an irreducible $D$ as
well, but this is not important for us.

Let $p=p_i$ be any of the points.
By \cite[II.1.5]{koll96} there is an algebraic space $U$ parametrizing
all morphisms
$\c\p^1\to X$ which pass through $p$ and have tangent space $\c\cdot v$
at
$p$. Let $F:U\times \p^1\to X$ be the universal morphism.
For $u\in U$ set $F_u:=F|_{\{u\}\times \p^1}$.

By \cite[IV.3.9]{koll96}
there is a $u\in U$ such that
$F_u^* T_X$ is
ample. (Strictly speaking, we need to apply the
refined version \cite[4.1.2]{lowdeg}.)
By composing $F_u$ with a multiple cover $\p^1\to \p^1$
we may assume that even $F_u^* T_X(-2)$ is
ample.
Hence there is a
connected open subset $W\subset U$ such that $F_w^* T_X(-2)$ is
ample for every $w\in W$. 

As in \cite[II.3.7]{koll96}
 this implies the following general position statements for $W$:
\begin{enumerate}
\item For any
$x\in X\setminus \{p\}$ the set 
$W_{x}:=\{w\in W | x\in \im F_w\}$
 has complex codimension $n-1$
in $W$.

\item For any
$x_1\neq x_2\in X\setminus \{p\}$ the set 
$W_{x_1,x_2}:=\{w\in W | x_1,x_2\in \im F_w\}$
 has complex codimension $2n-2$
in $W$.
\end{enumerate}

I claim that  $\im F_w$ satisfies
(\ref{gen.pos.lem}.1--2) for  general $w\in W$.  Indeed,
the set of all maps intersecting $Z$ is
$\cup_{x\in Z}W_x$ and this has codimension $\geq (n-1)-(n-2)=1$.
The set of those maps which intersect $X(\r)$ has real codimension
$\geq (2n-2)-n\geq 1$. 

Moreover, the 
set of those maps which pass through a conjugate point pair has 
codimension
$\geq (2n-2)-n\geq 1$, hence $\im (F_w)$ and $\overline{\im (F_w)}$
are disjoint for general $w\in W$. 

To complete the proof, first let $p=p_1$ and set
$D_1:=\im (F_w)\cup \overline{\im (F_w)}$
for  suitable general $w\in W$.

Next let $p=p_2$ and set $Z_1=Z\cup D_1$. As before we obtain $D_2$
and we can proceed to find the rest of the $D_i$ inductively.
Finally set $D=\cup_iD_i$.
\qed

\section{Torus bundles}

The aim of this section is to exhibit examples 
of real algebraic spaces  $X$ such that   $X$  is birational to $\p^3$
and $X(\r)$ is a torus bundle over a circle with prescribed monodromy.
We start with   some easy results about subgroups of $GL(2,\z)$.

\begin{say}[Subgroups of $GL(2,\z)$]{\ }\label{gl2.subgroups}

 Let $H_1\subset GL(2,\z)$ denote the subgroup of all matrices
$$
\left(
\begin{array}{cc}
a& b\\
c&d
\end{array}
\right)
\qtq{where $c\equiv 0\mod 2$.}
$$

 Let $H_2\subset GL(2,\z)$ denote the subgroup of all matrices
$$
\left(
\begin{array}{cc}
a& b\\
c&d
\end{array}
\right)
\equiv
\left(
\begin{array}{cc}
1& 0\\
0&1
\end{array}
\right)
\mod 2.
$$
The matrices 
$$
\left(
\begin{array}{cc}
0& 1\\
-1&-1
\end{array}
\right)
\qtq{and}
\left(
\begin{array}{cc}
0& 1\\
1&0
\end{array}
\right)
$$
generate a 6 element subgroup    $K_6\subset GL(2,\z)$ and
$GL(2,\z)$ is the semidirect product of $H_2$ by $K_6$.

Let $H_4\subset GL(2,\z)$ denote the subgroup of matrices
$$
\left(
\begin{array}{cc}
a& b\\
c&d
\end{array}
\right)
\qtq{where $a,d\equiv 1\mod 4$ and $b,c\equiv 0\mod 2$.}
$$
The matrices 
$
\left(
\begin{array}{cc}
\pm 1& 0\\
0&\pm 1
\end{array}
\right)
$
form a 4 element subgroup    $K_4\subset H_2$ and
$H_2$ is the semidirect product of  $H_4$ by $K_4$.
\end{say}

\begin{lem}\label{H4.gen.lem}
$H_4$ is generated by the matrices
$\left(
\begin{array}{cc}
1& 2\\
0&1
\end{array}
\right)
$
and
$
\left(
\begin{array}{cc}
1& 0\\
2&1
\end{array}
\right).
$
\end{lem}

Proof. Right multiplication by these matrices and their inverses changes
the first row
$(a,b)$ of a matrix to $(a,b\pm 2a)$ or to $(a\pm 2b,b)$. This allows
us to run a version of the Euclidean algorithm on the pair $(a,b)$. 
Taking into account that $a$ and $b$ are relatively prime, $a\equiv
1\mod 4$ and $b$ is even, we see that we can always reduce to the case
when the first row is $(1,0)$. Such a matrix in $H_4$ is automatically a
power of the second matrix above.\qed
\medskip

This easily implies the following:

\begin{cor}\label{H1.gen.cor}
$H_1$ is generated by 
$\left(
\begin{array}{cc}
\pm 1& *\\
0&\pm 1
\end{array}
\right)
$
and
$
\left(
\begin{array}{cc}
1& 0\\
2&1
\end{array}
\right).
$\qed
\end{cor}

\begin{say}[Birational transformations along ruled
surfaces]\label{ruled.bir.tr}{\ }

Let $X$ be a complex 3--fold and $F\subset X$
a ruled surface  with  
 normal bundle   $N_{F|X}$. Assume that $(E\cdot N_{F|X})=0$
where $E\subset F$ is a fiber of the projection $p:F\to C$.
Let $D\subset F$ be  a section of $p$; thus $D\cong C$ and
 $(D\cdot E)=1$.
Let $X'=B_DX$ be the blow up of $D$
with projection $g:X'\to X$. 
We blow up a curve inside $F$, so 
 the birational transform $g^{-1}_*(F)$ of $F$
is isomorphic to $F$.
If $N$ denotes the  normal bundle of $g^{-1}_*(F)\subset X'$
then $N\cong N_{F|X}(-D)$. Thus $(E\cdot N)=-1$. By (\ref{ruled.contr}) 
$g^{-1}_*(F)$ can be contracted and we obtain a birational transformation
$\rho:X\map Y$. $\rho$ is the composite of a blow up 
of the curve $D$ 
followed by the inverse of the blow up of a curve $C$.
As in (\ref{ruled.contr}) $Y$ need not be  projective, even if $X$ is.
\end{say}

\begin{say}[Birational transformations  of quadric bundles,
 local case]\label{bir.tr.quad.loc}{\ }

The smooth quadric surface $Q\cong \c\p^1\times \c\p^1$
is ruled in two distinct ways. Let $p:Q\to \c\p^1$ be the second
projection and   $E= \c\p^1\times (0:1)$  a fiber. Set  $F=
(0:1)\times
\c\p^1$  and let $D\subset Q$ denote a section of $p$.
$D$ is linearly equivalent to $mE+F$ for some $m\geq 0$.

Let $0\in \Delta$ be the complex unit disc
and consider the complex 3--manifold $X:=Q\times \Delta$
with second projection $\pi:X\to \Delta$.
Set $F=Q\times \{0\}$. Then $D\subset F\subset X$
satisfy the assumptions of (\ref{ruled.bir.tr}), hence we can construct
a birational transformation $\rho:X\map Y$.

 $(D\cdot D)=2m$
and so the normal bundle of $D$ in $X$ is $\o_D+\o_D(2m)$. Hence the
 exceptional divisor of the blow up 
 $g:X'\to X$ is isomorphic to
the ruled surface $F_{2m}$. We can now contract the birational transform
 $g^{-1}_*F$ to obtain  $Y$.
All these transformations take place over the central fiber of $\pi$,
so we still have a smooth, proper morphism
$\Pi:Y\to \Delta$. We have not changed the family 
over $\Delta\setminus\{0\}$ but the fiber of $\Pi$ over $0$ is now
$F_{2m}$.

Assume now that everything is defined over $\r$.
$X(\r)\to (-1,1)$ is a trivial torus bundle.
$Y(\r) \to (-1,1)$ is still smooth and proper, so it is again
a trivial torus bundle. We have, however, changed the trivialization.
 
In order to compute the change in the  trivialization, 
let $b^+,b^-\in (-1,1)$ be   a positive and a negative  point.
 Fix   bases
in $H_1(\pi^{-1}(b^{\pm}),\z)$ by  setting 
$$
e^{\pm}= E(\r)\times (0:1)\times \{b^{\pm}\}\qtq{and}
f^{\pm}=
(0:1)\times
F(\r)\times \{b^{\pm}\};
$$
 in both cases using the positive orientation on
$E(\r)$ and $F(\r)$. As we go from $b^-$ to $b^+$, parallel
transport produces an isomorphism
$$
T_X:H_1(\pi^{-1}(b^{-}),\z)\to H_1(\pi^{-1}(b^{+}),\z)
$$
whose marix is the identity.

$Y(\r)\to (-1,1)$ is another torus bundle, so we get another
parallel
transport isomorphism
$$
T_Y:H_1(\Pi^{-1}(b^{-}),\z)\to H_1(\Pi^{-1}(b^{+}),\z).
$$
Let us compute the 
matrix of $T_Y$. 

The homology class of $D(\r)\times \{b^{\pm}\}$
 is of the form  $ne^{\pm}+f^{\pm}$ for some $|n|\leq m$.
It is easy to see that any $n$ can occur for suitable choice of $D$.
For our purposes we need only the cases $n\in \{-1,0,1\}$.
These are realized by the skew diagonal, $F$ and the diagonal.

$D\times (-1,1)$ intersects $F$ transversally in $D\times \{0\}$,
so $g^{-1}_*(D\times (-1,1))\cong D\times (-1,1)$ is disjoint from
$g^{-1}_*(Q)$. Thus  $T_Y$ transports $D(\r)\times \{b^-\}$ to
$D(\r)\times \{b^+\}$  and preserves the orientation,  so
$$
T_Y
\left(
\begin{array}{c}
n\\
1
\end{array}
\right)=
\left(
\begin{array}{c}
n\\
1
\end{array}
\right).
$$
 Consider now $E\times (-1,1)$. We
blow up the unique point $(E\cap D)\times \{0\}$ and then contract the
birational transform of $E\times \{0\}$. 
Let  $t$ be a local coordinate on $\Delta$ and $(x_0:x_1)$
   coordinates on $E$ such that $E\cap D=(0:1)$. Then the curves
$x_0=\lambda x_1$ give the parallel transport before blow up
and the curves
$x_0=\lambda tx_1$ give the parallel transport after blow up.
Thus 
 $T_Y$ transports $E(\r)\times \{b^-\}$ to
$E(\r)\times \{b^+\}$  but it  reverses the orientation.
Therefore
$$
T_Y
\left(
\begin{array}{c}
1\\
0
\end{array}
\right)=
\left(
\begin{array}{c}
-1\\
0
\end{array}
\right)
\qtq{and so} 
T_Y=
\left(
\begin{array}{cc}
-1& 2n\\
0&1
\end{array}
\right).
$$
Of course we could have used the other 
projection of the quadric as a ruled surface structure. These give
parallel transport matrices
$$
T'_Y=
\left(
\begin{array}{cc}
1& 0\\
2n&-1
\end{array}
\right).
$$
\end{say}

\begin{say}[Birational transformations  of quadric bundles, global
case]{\ }

Let $X=Q\times \c\p^1$ with projection $\pi:X\to \c\p^1$.
We can pick  points $p_i\in \r\p^1$
and at each of these points perform one of the  birational
transformations discussed above. At the end we obtain
a   real algebraic space $Y$
such that the projection $\Pi:Y(\r)\to \r\p^1$
is a torus bundle. 

The  monodromy action on $H_1$ can be written as the product of
the local parallel transport matrices $T_Y$ and $T'_Y$. 
Observe that
$$
\left(
\begin{array}{cc}
-1& \pm 2\\
0&1
\end{array}
\right)
\left(
\begin{array}{cc}
-1& 0\\
0&1
\end{array}
\right)
=
\left(
\begin{array}{cc}
1& \pm 2\\
0&1
\end{array}
\right).
$$
and
$$
\left(
\begin{array}{cc}
1& 0\\
0&-1
\end{array}
\right)
\left(
\begin{array}{cc}
1& 0\\
\pm 2&-1
\end{array}
\right)
=
\left(
\begin{array}{cc}
1& 0\\
\pm 2&1
\end{array}
\right).
$$
Thus by (\ref{gl2.subgroups}) and (\ref{H4.gen.lem})
the matrices
$$
\left(
\begin{array}{cc}
-1& 0\\
0&1
\end{array}
\right), 
\left(
\begin{array}{cc}
-1& \pm 2\\
0&1
\end{array}
\right),
\left(
\begin{array}{cc}
1& 0\\
0&-1
\end{array}
\right),
\left(
\begin{array}{cc}
1& 0\\
\pm 2&1
\end{array}
\right)
$$
generate $H_2\subset GL(2,\z)$.
Hence
we   get    examples of real algebraic spaces  $Y$
such that $Y$ is birational to $\p^3$ and $Y(\r)$ is a torus bundle over
a circle with arbitrary  monodromy in $H_2$.
\end{say}

Unfortunately, the  example we need in (\ref{ben-mar.B.thm}.2)
can not be realized with monodromy in $H_2$. Thus we need
to introduce a further twist in the construction.

\begin{say}[An example of a degree 6 Del Pezzo surface]{\
}\label{d6DP.defn}

Let $Q$ be the quadric surface $(x^2+y^2+z^2=t^2)\subset \p^3$
and $p_1,p_2,p_3\in Q(\r)$  real points.
 Let $C\subset Q$
denote the intersection of $Q$ with the plane spanned by the
points $p_1,p_2,p_3$.  Let  us blow up the 3 points
to obtain $h:T\to Q$
with  exceptional curves $E_i$. $h^{-1}_*(C)$ is a smooth rational
curve with selfintersection $-1$, thus it can be  contracted
$g:T\to  S$. $S$ is a degree 6 Del Pezzo surface.
Set $F_i:=g_*(E_i)$. 

$T(\r)$ is the connected sum of 3 copies of
$\r\p^2$. The birational transform of $C(\r)$ intersects
 each $E_i(\r)$ in a single point, thus $S(\r)$ is a torus.

Fix an orientation on
$Q(\r)\sim S^2$. This induces an orientation on each
 $E_i(\r)$.
$C(\r)$ is a  circle which divides $Q(\r)$ into two discs.
The image of either of these discs in $S(\r)$ has
$\cup F_i(\r)$ as its boundary. Thus $\sum [F_i(\r)]=0\in
H_1(S(\r),\z)$. The intersection number of $F_i(\r)$ and
of $F_j(\r)$ is $\pm 1$, hence any two of the 
$[F_i(\r)]$ form a basis of $H_1(S(\r),\z)$.

Assume now that we have a diffeomorphism $Q(\r)\to Q(\r)$
which preserves orientation and permutes the three points 
$(p_1,p_2,p_3)\mapsto (p_2,p_3,p_1)$.
Using  $[F_1(\r)], -[F_2(\r)]$
as a basis, 
the induced map on $H_1(S(\r),\z)$ is  given by  the
matrix
$$
\tau_3=\  \left(
\begin{array}{cc}
0& 1\\
-1&-1
\end{array}
\right).
$$
\end{say}

\begin{lem}\label{K6.monodr.lem}
 Let $K_6\subset GL(2,\z)$ be the subgroup defined in
(\ref{gl2.subgroups}). For every $M\in K_6$ there is a smooth,
projective, real algebraic variety $X$ such that
\begin{enumerate}
\item there is a morphism $X\to \p^1$ which is an analytically
locally trivial fiber bundle  over 
$\r\p^1$,
\item $X(\r)\to \r\p^1$ is a torus bundle with monodromy $M$, and
\item $X$ is obtained from $\p^2\times \p^1$ by a sequence of 
smooth blows ups and  blow downs.
\end{enumerate}
\end{lem}

Proof.  First we give an example with order 3 monodromy
in $K_6$. 
Fix  a smooth conic $C\subset  Q$ and consider
the surface $C\times \p^1$. Let $D\subset C\times \p^1$
be the graph
of the morphism
$C\cong \p^1 {\to} \p^1$
given as $z\mapsto (z^3-4z)/(z^2-1)$. 
The induced projection $D\to \p^1$
has degree 3 and  $D(\r)\to \r\p^1$
is three--to--one over every point.
Let $Y\to Q\times \p^1$ be the blow up of
$D$. The birational transform of $C\times \p^1$
is contractible in $Y$. By contracting it we obtain
$Y\to X\to \p^1$.

Over each real point we have an example of the construction
of a degree 6 Del Pezzo surface given in (\ref{d6DP.defn}),
thus $X(\r)\to \r\p^1$ is a torus bundle with
 order 3 monodromy $\tau_3$ as above.
The quadric $Q$ is obtained from $\p^2$ by two blow ups and
one blow down,
so $Q\times \p^1$  is obtained from $\p^2\times \p^1$ by two smooth blow ups
and one smooth blow down.

Examples with  order 2 monodromy can also be obtained from
degree 6 Del Pezzo surfaces, but they are easier to get from a
family of  quadrics. We again start with $C\times \p^1$
but we realize $C$ as a line in $\p^2$. Choose
 $D\subset C\times \p^1$ as the graph
of the morphism
$C\cong \p^1 {\to} \p^1$
given as $z\mapsto (z^2-1)/z$. 
  The previous blow up and
blowdown construction gives an example with monodromy
$$
  \left(
\begin{array}{cc}
0& 1\\
-1& 0
\end{array}
\right).
$$
We have to perform one more transformation as in
(\ref{bir.tr.quad.loc}) with $n=0$ to get an example with
monodromy
$$
\tau_2=\   \left(
\begin{array}{cc}
0& 1\\
1& 0
\end{array}
\right).
$$
The    matrices $\tau_3$ and $\tau_2$, together with
their conjugates, give all elements of $K_6$.\qed

\begin{say}[Birational transformations along degree 6 Del Pezzo
surfaces]\label{deg6DP.bir.tr}{\ }

Let $X$ be a complex 3--fold with real structure and $F\subset X$
a   surface  isomorphic to $S$ defined in (\ref{d6DP.defn}). 
  Assume that  $F$ has trivial normal bundle   $N_{F|X}\cong \o_F$.
Then $F(\r)$ looks like a ruled surface. We would like to
perform the same transformations we did in (\ref{ruled.bir.tr}).
In order to do this, we need to perform some preparatory
transformations first. During these we do not change the
real points but the surface $F$ becomes a ruled surface with the
requisite normal bundle. This is done in two steps.

Through each point $p_i\in Q$ there is a conjugate pair of lines.
Let $G_i\subset F\cong S$ denote their birational transform.

First blow up $G_1\subset F$ to get $X_1$ and let $F_1\subset X_1$ be
the
birational transform of $F$. 
Then $F_1\cong F$ and $G_2\subset
F_1\subset X_1$ 
is a conjugate pair of smooth rational curves
with normal bundles $\o_{\p^1}(-1)+\o_{\p^1}(-1)$.  Thus we can 
 flop $G_2\subset X_1$ to obtain
$X_2$.  Let $F_2\subset X_2$ denote the birational transform of $F_1$.
$F_2$ is obtained from $F_1\cong S$ by contracting $G_2$,
thus $F_2$ can be identified with $\p^1\times \p^1$.
The transformations $X\map X_1\map X_2$ are isomorphisms
along the real points. 
Hence we are in the situation of
(\ref{ruled.bir.tr}) and we can perform the transformations
discussed there. 
\end{say}

Putting all of these together, we obtain the main result of the section:

\begin{thm}\label{GL2.monodr.thm}
 For every $M\in GL(2,\z)$ there is a smooth,
 real, algebraic space $X$ such that
\begin{enumerate}
\item $X(\r)$ is a torus bundle over $S^1$ with monodromy $M$, and
\item $X$ is obtained from  $\p^2\times \p^1$ by a sequence of
smooth blow ups and blow downs.
\end{enumerate}
\end{thm}

Proof. Write $M=M_1M_2$ where $M_1\in K_6$ and $M_2\in H_2$. 
We also write $M_2=\prod T_i$ where each $T_i$ is one of the matrices
$$
\left(
\begin{array}{cc}
-1& 0\\
0&1
\end{array}
\right), 
\left(
\begin{array}{cc}
-1& \pm 2\\
0&1
\end{array}
\right),
\left(
\begin{array}{cc}
1& 0\\
0&-1
\end{array}
\right),
\left(
\begin{array}{cc}
1& 0\\
\pm 2&1
\end{array}
\right).
$$
By (\ref{K6.monodr.lem}) we can obtain an $X_1$ as in 
(\ref{GL2.monodr.thm}) with monodromy $M_1$. 
Pick points in cyclic order $p_i\in \r\p^1$.
$X_1\to \p^1$ is a trivial $S$-bundle over an analytic  neighborhood of
$p_i$, thus we can perform birational transformations
as in (\ref{deg6DP.bir.tr}) to change the global monodromy by the matrix
$T_i$. Doing this at all points $p_i$, we obtain $X\to \p^1$.
The   monodromy of $X(\r)\to \r\p^1$ is thus 
$M_1\prod T_i=M$. 
\qed

\begin{rem} From (\ref{GL2.monodr.thm}) we only need the special case 
when
$$
M=
\left(
\begin{array}{cc}
0& 1\\
1&1
\end{array}
\right)=
\left(
\begin{array}{cc}
0& 1\\
-1&-1
\end{array}
\right)
\left(
\begin{array}{cc}
-1& -2\\
0&1
\end{array}
\right).
$$
In this case we need to perform the birational tranformations in 
a single fiber only. It is not hard to see that during this process we
create an effective curve which is homologous to zero,
so the resulting algebraic space $X$ is  not projective.
\end{rem}


\noindent Princeton University, Princeton NJ 08544-1000

\begin{verbatim}kollar@math.princeton.edu\end{verbatim}

\end{document}